\documentclass[a4paper]{article}

\usepackage{amsmath,amsfonts,amssymb}
\usepackage{amsthm}

\newcommand{\mathsym}[1]{{}}

\theoremstyle{plain}
\newtheorem{theorem}{Theorem}
\newtheorem{proposition}[theorem]{Proposition}
\newtheorem{lemma}[theorem]{Lemma}

\newtheorem{remark}{Remark}[section]

\begin{document}
\title{Maximal $m$-distance sets containing the representation of the Johnson graph $J(n, m)$}
\author{Eiichi Bannai${}^{*}$, Takahiro Sato${}^{\dagger}$, Junichi Shigezumi${}^{\ddagger}$}
\date{}

\maketitle \vspace{-0.1in}
\begin{center}

${}^{*}$Department of Mathematics, Shanghai Jiao Tong University\\
800 Dongchuan Road, Shanghai  200240, China\\
{\it E-mail address} : bannai@sjtu.edu.cn\vspace{0.05in}\\

${}^{\dagger}$Oita Maizuru High School\\
1-19-1 Imazuru, Oita 870-0938, Japan\vspace{0.05in}\\

${}^{\ddagger}$Institute of Mathematics for Industry, Kyushu University\\
744 Motooka, Nishi-ku, Fukuoka 819-0395, Japan\\
{\it E-mail address} : j1.shigezumi@gmail.com \vspace{-0.05in}
\end{center} \quad

\begin{quote}
{\small\bfseries Abstract.}
We classify the maximal $m$-distance sets in $\mathbb{R}^{n - 1}$ which contain the representation of the Johnson graph $J(n, m)$ for $m = 2, 3$.
Furthermore, we determine the necessary and sufficient condition for $n$ and $m$
such that the representation of the Johnson graph $J(n, m)$ is not maximal as an $m$-distance set.

Also, we classify the maximal two-distance sets in $\mathbb{R}^{n - 1}$ which contain the representation of $J(n - 1, 2)$.\\  \vspace{-0.15in}

 \vspace{0.15in}\\

\end{quote}

\section*{Introduction}

The purpose of our study is to find an `interesting finite subset' of the Euclidean space.
$s$-{\it distance sets}, where $s$ is the number of distances between distinct vectors,
sometimes provide such an interesting subset.
In particular, many interesting examples exist of the {\it maximum} $s$-distance sets,
those having the largest cardinality of $s$-distance sets,
including the regular pentagon and heptagon in $\mathbb{R}^2$,
and the regular octahedron and icosahedron in $\mathbb{R}^3$.

However, determining the maximum $s$-distance sets is not easy in general.
For $s$-distance sets in $\mathbb{R}^2$, we have a classification on at most five-distance sets
(by Kelly\cite{K}, Erd\"os-Fishburn\cite{E-F} and Shinohara\cite{Sh1}).
For $\mathbb{R}^3$, we have a classification on at most three-distance sets
(by Croft\cite{C}, Einhorn-Schoenberg\cite{E-S2} and Shinohara\cite{Sh2}, \cite{Sh3}).

On the other hand, for the maximum two-distance sets, we have the following table: 
\begin{center}
\begin{tabular}{c|cccccccc}
\hline
dimension & $1$ & $2$ & $3$ & $4$ & $5$ & $6$ & $7$ & $8$\\
\hline
\hline
max. card. & $3$ & $5$ & $6$ & $10$ & $16$ & $27$ & $29$ & $45$\\
\hline
\# of sets & $1$ & $1$ & $6$ & $1$ & $1$ & $1$ & $1$ & $\geqslant 1$\\
\hline
\end{tabular}
\end{center}
For the dimensions $4 \leqslant n \leqslant 6$, Seidel\cite{Se} conjectured the values above.
In 1997, Lison\v{e}k\cite{L} proved the maximum cardinalities of two-distance sets in $\mathbb{R}^n$ for $4 \leqslant n \leqslant 8$.
Furthermore, he found examples for $n = 7, 8$, and proved the uniqueness of the maximum two-distance set for $4 \leqslant n \leqslant 7$.

Here, for positive integers $n$ and $m$, the {\it Johnson graph} $J(n, m)$ has vertex set $V(n, m)$ and edge set $E(n, m)$ as follows:
\begin{align*}
V(n, m) &:= \left\{  \left\{ i_1, \ldots, i_m \right\} : 1 \leqslant i_1 < \cdots < i_m \leqslant n \right\},\\
E(n, m) &:= \left\{ (v_i, v_j) : | v_i \cap v_j | = m, \; v_i, v_j \in V(n, m) \right\}.
\end{align*}
We note that the two-distance sets in $\mathbb{R}^n$ for $n = 7, 8$ contain
the representations of the Johnson graphs $J(7, 2)$ and $J(9, 2)$, respectively.
Both of these sets are {\it maximal} two-distance sets
which contain the representations of the corresponding Johnson graphs
in the each Euclidean space $\mathbb{R}^n$.
Here, an $s$-distance set is said to be {\it maximal}
if no other $s$-distance set contains it in the given space.

We have the following representation of the Johnson graph $J(n, m)$ in Euclidean space $\mathbb{R}^{n - 1}$:
\begin{equation}
\tilde{J}(n, m) = (1^m, 0^{n - m})^P.
\end{equation}
Here, the exponents inside the parentheses indicate the number of occurrences of the corresponding numbers,
and the exponent $P$ outside indicates that we should take every permutation.
We note that every vector should be on the hyperplane $\{ (x_1, \ldots, x_n) \in \mathbb{R}^n :  \sum_i x_i = m \}$.
In addition, since we have $\tilde{J}(n, m) \simeq \tilde{J}(n, n - m)$, we may assume $n \geqslant 2 m$.
Then, $\tilde{J}(n, m)$ is an $m$-distance set.\\

In this paper, we investigate the maximal $m$-distance sets which contain the representation $\tilde{J}(n, m)$ of the Johnson graph $J(n, m)$ in $\mathbb{R}^{n}$ and $\mathbb{R}^{n-1}$.

The first problem is to determine maximal $m$-distance sets in $\mathbb{R}^{n - 1}$ which contain $\tilde{J}(n, m)$,
which correspond to the maximum two-distance sets in $\mathbb{R}^8$ with $J(9, 2)$ introduced by Lison\v{e}k.
We classify the maximal $m$-distance sets which contain $\tilde{J}(n, m)$ for $m = 2, 3$.
Furthermore, we determine the necessary and sufficient condition for $n$ and $m$,
where $\tilde{J}(n, m)$ is not maximal as an $m$-distance set.

We have the following theorem.

\begin{theorem}\label{th-main}
The representation $\tilde{J}(n, m)$ of the Johnson graph $J(n, m)$ in the Euclidean space $\mathbb{R}^{n - 1}$
is not maximal as an $m$-distance set
if and only if  it satisfies the following conditions:
\begin{equation}
n > n_0 \qquad \text{and} \qquad 3 n_0 - \frac{n_0^2}{n} \leqslant 4 m. \label{cond-1}
\end{equation}
Furthermore, if the pair $(n, m)$ satisfies the above conditions,
we can add each vector of the following sets to $\tilde{J}(n, m)$ while maintaining $m$-distance:
\begin{equation}
\begin{cases}
\left( \left( \frac{n_0}{n} \right)^{n - n_0 + m}, \; \left( - 1 + \frac{n_0}{n} \right)^{n_0 - m} \right)^P & m < n_0,\\
\qquad \left( \left( \frac{m}{n} \right)^n \right) & m = n_0,\\
\left( \left( 1 + \frac{n_0}{n} \right)^{m - n_0}, \; \left( \frac{n_0}{n} \right)^{n + n_0 - m} \right)^P & m > n_0.
\end{cases}
\end{equation}
\end{theorem}

Here, $n_0$ is the special factor of $n$. First, we consider the factorization of integer $n$:
\begin{equation}
n = 2^{e_0} \prod_{i > 0} p_i^{e_i}, \quad (p_i: \, \text{odd prime}, \; e_i \in \mathbb{Z}).
\end{equation}
Then, we have the integer $n_0$ as follows:
\begin{equation}
n_0 :=
\begin{cases}
\qquad {\displaystyle \prod_{i > 0} p_i^{\lceil e_i / 2 \rceil}} & e_0 = 0,\\
{\displaystyle 2^{\lceil (e_0 + 1) / 2 \rceil}  \prod_{i > 0} p_i^{\lceil e_i / 2 \rceil}} & e_0 > 0.
\end{cases}
\end{equation}

The second problem is to determine maximal $m$-distance sets in $\mathbb{R}^{n - 1}$ which contain $\tilde{J}(n - 1, m)$,
which corresponds to the maximum two-distance sets in $\mathbb{R}^7$ with $J(7, 2)$ introduced by Lison\v{e}k.
We classify the maximal two-distance sets which contain $\tilde{J}(n, 2)$.

In Section \ref{sec-m_2_3}, we classify the maximal $m$-distance sets in $\mathbb{R}^{n - 1}$ which contain $\tilde{J}(n, m)$ for $m = 2, 3$.
In Section \ref{sec-pf-th-main}, we give the proof of Theorem \ref{th-main}.
We show some general results related to Theorem \ref{th-main} in Section \ref{sec-rem}.
Finally, in Section \ref{sec-n-1}, we classify the maximal two-distance sets in $\mathbb{R}^{n - 1}$ which contain $\tilde{J}(n - 1, 2)$.

\section{Preliminaries}\label{sec-pre}

In the definition of $\tilde{J}(n, m)$, we take every permutation.
Thus, if we can add a vector $x = (x_1, \ldots, x_n)$ to $\tilde{J}(n, m)$ while maintaining $m$-distance,
then we can also add each vector of the set $(x_1, \ldots, x_n)^P$ to $\tilde{J}(n, m)$.
Furthermore, we have
\begin{equation*}
\left\{ d(x, y) : x, y \in \tilde{J}(n, m) \right\} = \left\{ \sqrt{2 i} : i = 1, 2, \ldots, m \right\},
\end{equation*}
where $d(x, y)$ is the Euclidean distance between $x$ and $y$.
Then it is easy to show that,
for any vector $x = (x_1, \ldots, x_n)$ which we can add to $\tilde{J}(n, m)$ while maintaining $m$-distance,
the number of the elements of $x$ should be at most $m$,
and the differences between the elements of $x$ should be integers which are at most $m - 1$.
In conclusion, we have the following notation for such a vector:
\begin{equation}
X = \left( \left( 1 - \frac{k_0}{n} \right)^{k_1}, \; \left( - \frac{k_0}{n} \right)^{k_2}, \; \left( -1 - \frac{k_0}{n} \right)^{k_3}, \; \ldots, \; \left( - (l - 2) - \frac{k_0}{n} \right)^{k_l}  \right)^P. \label{not-adpt}
\end{equation}
Here, we have
\begin{gather}
\notag l \leqslant m, \qquad
\forall j > 0 \quad k_j \geqslant 0,\\
\sum_{j>0} k_j = n, \qquad \sum_{j>0} j \, k_j = 2 n - k_0 - m. \label{cond-adpt}
\end{gather}
The final condition is that such vectors must be on the hyperplane which contains $\tilde{J}(n, m)$.

Now, we take vectors $x = (x_1, \ldots, x_n) \in \tilde{J}(n, m)$ and $y = (y_1, \ldots, y_n) \in X$.
We give the following definition.
\begin{equation*}
i_{j} := \left| \left\{ 1 \leqslant k \leqslant n \; ; \; x_k = 1, \; y_k = - (j - 2) - \frac{k_0}{n} \right\} \right|.
\end{equation*}
The distance $d(x, y)$ between the vectors $x$ and $y$ depends on only the combinations of the elements $x_k$ and $y_k$.
Thus, we have
\begin{align}
\notag \left( d(x, y) \right)^2 &= \sum_j i_j \, \left\{ 1 - \left( - (j - 2) - \frac{k_0}{n} \right)\right\}^2
   + \sum_j (k_j - i_j) \, \left\{ - \left( - (j - 2) - \frac{k_0}{n} \right)\right\}^2\\
 &= 4 k_0 + 3 m - 4 n - \frac{k_0^2}{n} + \sum_j j^2 \, k_j + 2 \sum_j (j - 1) \, i_j. \label{eq-dis^2}
\end{align}
Note that, for another element $x' \in X$, the difference between $\left( d(x, y) \right)^2$ and $\left( d(x', y) \right)^2$
is just the difference in the final term $2 \sum_j (j - 1) \, i_j$ in the above equation.

\section{Cases $m=2, 3$}\label{sec-m_2_3}

\subsection{Case $m = 2$}
We may assume $n \geqslant 2 m = 4$.
As we showed in the previous section,
if we have some vector which we can add $\tilde{J}(n, 2)$ while maintaining two-distance,
then we have the following set which contains the vector:
\begin{equation*}
X = \left( \left( 1 - \frac{k_0}{n} \right)^{k_1}, \; \left( - \frac{k_0}{n} \right)^{k_2} \right)^P.
\end{equation*}
By condition (\ref{cond-adpt}), we have $k_1 = k_0 + 2$, $k_2 = n - k_0 - 2$.
Then, for the vectors $x \in \tilde{J}(n, 2)$ and $y \in X$,
by  (\ref{eq-dis^2}), we have
\begin{equation*}
\left( d(x, y) \right)^2 = k_0 - \frac{k_0^2}{n} + 2 \, i_2.
\end{equation*}

Now, we may assume that $-2 \leqslant k_0 < n - 2$
since $k_1$ and $k_2$ are nonnegative
and we have $X = \left( \left( \frac{2}{n} \right)^n \right)$ for both cases $k_0 = -2$ and $k_0 = n - 2$.
Then, we consider the following separate cases:
$k_0 = -2$, $k_0 = -1$, $0 \leqslant k_0 \leqslant n - 4$, and $k_0 = n - 3$.

For the cases $k_0 = -2$ and $k_0 = -1$, the squared distance $\left( d(x, y) \right)^2$ is not an even integer for any $n \geqslant 4$.
For the case $0 \leqslant k_0 \leqslant n - 4$, $i_2$ takes three values: $0$, $1$, and $2$.
Then, the number of the distances is three.
This contradicts the two-distance assumption.

Finally, for the case $k_0 = n - 3$, $i_2$ takes two values, $0$ and $1$,
and $\left( d(x, y) \right)^2 = n - 3 - \frac{(n - 3)^2}{n}$ and $n - 1 - \frac{(n - 3)^2}{n}$.
Since the two distances must be even integers at most $4$, we have only the case $n = 9$:
\begin{equation*}
X = \left( \left( \frac{1}{3} \right)^8, \; - \frac{2}{3} \right)^P.
\end{equation*}
In addition, the distances between any distinct vectors of $X$ are $\sqrt{2}$ and $2$.
Thus, we can add all the vectors of $X$ to $\tilde{J}(9, 2)$ while maintaining two-distance.

In conclusion, if $n \ne 9$, then $\tilde{J}(n, 2)$ is maximal as a two-distance set.
On the other hand, if $n = 9$, the maximal two-distance set which contains $\tilde{J}(9, 2)$ is the following:
\begin{equation}
\tilde{J}(9, 2) \cup X = \left( 1^2, \; 0^7 \right)^P \, \bigcup \, \left( \left( \frac{1}{3} \right)^8, \; - \frac{2}{3} \right)^P.
\end{equation}
This set is a maximum two-distance set with $45$ vectors in $\mathbb{R}^8$, which is the same as the set introduced by Lison\v{e}k.

\subsection{Case $m = 3$}
We may assume $n \geqslant 6$. We have
\begin{equation*}
X = \left( \left( 1 - \frac{k_0}{n} \right)^{k_1}, \; \left( - \frac{k_0}{n} \right)^{k_2}, \; \left( -1 - \frac{k_0}{n} \right)^{k_3} \right)^P.
\end{equation*}
Also, we have $k_2 = n + k_0 + 3 - 2 k_1$, $k_3 = - k_0 - 3 + k_1$.
Then, for the vectors $x \in \tilde{J}(n, 3)$ and $y \in X$, we have
\begin{equation*}
\left( d(x, y) \right)^2 = - k_0 + 2 k_1 - 6 - \frac{k_0^2}{n} + 2 i_2 + 4 i_3.
\end{equation*}

Here, we split the problem into the following four cases:
$(a)$ $k_2 = k_3 = 0$, $(b)$ $k_1, k_2 > 0$, $k_3 = 0$, $(c)$ $k_1, k_3 > 0$, $k_2 = 0$, and $(d)$ $k_1, k_2, k_3 > 0$.

For the case $(a)$, we have $(i_1, i_2, i_3) = (3, 0, 0)$. Then, we have the only case ``$k_0 = 6$, $k_1 = 9$, $n = 9$.''

For the case $(b)$, we consider the following separate cases:
$k_0 = 1$, $k_0 = 2$, $3 \leqslant k_0 \leqslant n-3$, $k_0 = n - 2$, and $k_0 = n-1$.
Similarly to the case $m = 2$, we can limit the case to ``$k_0 = 4$, $k_1 = 7$, $n = 8$.''

For the case $(c)$, we have $k_1 = \frac{n + k_0}{2}$, $k_3 = \frac{n - k_0}{2}$.
We split the problem into several cases in terms of $k_3$.
However, we have no case which satisfies the conditions $n \geqslant 6$ and $\left( d(x, y) \right)^2$ are even integers.

Finally, for the case $(d)$, we note the fact that the difference between
$\max_{X} \left( d(x, y) \right)^2$ and $\min_{X} \left( d(x, y) \right)^2$ is at most $4$.
Thus, we also have $\max_{i_2, i_3} (2 i_2 + 4 i_3) - \min_{i_2, i_3} (2 i_2 + 4 i_3) \leqslant 4$.
We have the case $(i_1, i_2, i_3) = (1, 1, 1)$, so we have $0 < 2 i_2 + 4 i_3 < 12$.
Then, we have $k_1 \leqslant 2$, $k_3 \leqslant 2$.

We split the problem into several cases in terms of $k_3$ and $k_2$,
and we check the number $2 i_2 + 4 i_3$ for all the possible triples $(i_1, i_2, i_3)$.
Since we have $n \geqslant 6$ and the fact that $\left( d(x, y) \right)^2$ must be even integers at most $6$,
we can limit the case to ``$k_0 = -3$, $k_1 = 1$, $n = 9$.''

In conclusion, we have the only three cases, namely,
$(i)$ $k_0 = 4$, $k_1 = 7$, $n = 8$; $(ii)$ $k_0 = 6$, $k_1 = 9$, $n = 9$; and $(iii)$ $k_0 = -3$, $k_1 = 1$, $n = 9$.
Then, the corresponding sets $X$ are the following:
\begin{equation*}
X^{(i)} = \left( \left( \frac{1}{2} \right)^7, -\frac{1}{2} \right)^P, \quad
X^{(ii)} = \left( \left( \frac{1}{3} \right)^9 \right), \quad
X^{(iii)} = \left( \frac{4}{3}, \left( \frac{1}{3} \right)^7, -\frac{2}{3} \right)^P.
\end{equation*}
In addition, the distances between any distinct vectors of $X^{(i)}$ are $\sqrt{2}$,
and distances between the vector of $X^{(ii)}$ and any vectors of $X^{(iii)}$ are $\sqrt{2}$.
However, for $x_0 \in X^{(ii)}$ and every vector $y \in X^{(iii)}$, we have $y' = 2 x_0 - y \in X^{(iii)}$ and $d(y, y') = 2 \sqrt{2}$.
We can classify all the vectors of $X^{(iii)}$ into such pairs.
Then, we have the following set:
\begin{equation*}
{X^{(iii)}}' = \left\{ (x_1, \ldots, x_9) \in X^{(iii)} : x_i = \frac{4}{3}, x_j = -\frac{2}{3}, i > j \right\}.
\end{equation*}
This is one of the maximal subsets of $X^{(iii)}$ as a three-distance set.

In conclusion, if $n \ne 8, 9$, then $\tilde{J}(n, 3)$ are maximal as three-distance sets.
On the other hand, if $n = 8$, the maximal three-distance set which contains $\tilde{J}(8, 3)$ is the following:
\begin{equation}
\tilde{J}(8, 3) \cup X^{(i)},
\end{equation}
which has $64$ vectors.

For the case $n = 9$, we can choose the vector of $X^{(ii)}$
and half of the vectors of $X^{(iii)}$ by taking one vector from each pair $(y, y')$.
One of the maximal three-distance sets which contains $\tilde{J}(9, 3)$ is the following:
\begin{equation}
\tilde{J}(9, 3) \cup X^{(ii)} \cup {X^{(iii)}}'.
\end{equation}
Each of the maximal three-distance sets has $121$ vectors.

\section{Proof of Theorem \ref{th-main}}\label{sec-pf-th-main}

By the same reasoning as in the previous sections,
if we have some vectors which we can add to $\tilde{J}(n, m)$ while maintaining $m$-distance,
we have the following set which contains the vector:
\begin{equation*}
X = \left( \left( 1 - \frac{k_0}{n} \right)^{k_1}, \; \left( - \frac{k_0}{n} \right)^{k_2}, \; \left( -1 - \frac{k_0}{n} \right)^{k_3}, \; 
\ldots, \; \left( - (l - 2) - \frac{k_0}{n} \right)^{k_l}  \right)^P.
\end{equation*}
Furthermore, for the vectors $x \in \tilde{J}(n, m)$ and $y \in X$,
we have
\begin{equation*}
\left( d(x, y) \right)^2 = 4 k_0 + 3 m - 4 n - \frac{k_0^2}{n} + \sum_j j^2 \, k_j + 2 \sum_j (j - 1) \, i_j,
\end{equation*}
where $i_{j} = \left| \left\{ 1 \leqslant k \leqslant n \; ; \; x_k = 1, \; y_k = - (j - 2) - \frac{k_0}{n} \right\} \right|$. 

Let $y$ be a vector of $X$. Then we denote the maximum distance of vectors from $\tilde{J}(n, m)$ and $X$ as follows:
\begin{equation*}
M_X := \max_{\begin{subarray}{l}x \in \tilde{J}(n, m)\end{subarray} } \left( d(x, y) \right)^2.
\end{equation*}
We note that this does not depend on the choice of the vector $y \in X$.

Recall that, for another choice, elements $x' \in X$ and $y' \in \tilde{J}(n, m)$,
the difference between square distances $\left( d(x, y) \right)^2$ and $\left( d(x', y') \right)^2$
appears in the final term, $2 \sum_j (j - 1) \, i_j$, of the above equation, which is an even integer.
Thus, if the maximum $M_X$ is an even integer at most $2 m$,
then $\left( d(x, y) \right)^2$ is also an even integer at most $2 m$ for every $x \in X$ and $y \in \tilde{J}(n, m)$.
Now, we have the following lemma.

\begin{lemma} \label{lem-mx}
We have some vectors $x$ which we can add to $\tilde{J}(n, m)$ while maintaining $m$-distance
if and only if $M_X$ is an even integer at most $2 m$ for the corresponding set X, with $x \in X$.
\end{lemma}

In addition, we can easily show that we have some $1 \leqslant j_0 \leqslant l$ such that
\begin{equation*}
i_j =
\begin{cases}
0 & j < j_0,\\
m - \sum_{j > j_0} k_j & j = j_0,\\
k_j & j > j_0.
\end{cases}
\end{equation*}
Furthermore, we have some vectors $x \in X$ and $y \in \tilde{J}(n, m)$ which correspond to the above $i_j$.
This tuple $\{ i_j \}_j$ maximizes $2 \sum_j (j - 1) \, i_j$, and thus it also maximizes $\left( d(x, y) \right)^2$.
We define $I_X := 2 \sum_j (j - 1) \, i_j$ for the above tuple $\{ i_j \}_j$, which corresponds to $M_X$.

We consider another set $X'$ from the set $X$:
\begin{equation*}
X' := \left( \left( 1 - \frac{k_0}{n} \right)^{{k_1}'}, \; \left( - \frac{k_0}{n} \right)^{{k_2}'}, \; \left( -1 - \frac{k_0}{n} \right)^{{k_3}'}, \; 
\ldots, \; \left( - (l - 2) - \frac{k_0}{n} \right)^{{k_l}'}  \right)^P.
\end{equation*}
where, if $l > 3$, then ${k_1}' = k_1 - 1$, ${k_2}' = k_2 + 1$, ${k_{l - 1}}' = k_{l - 1} + 1$, ${k_l}' = k_l - 1$,
and ${k_j}' = k_j$ for $3 \leqslant j \leqslant l - 2$.
On the other hand, if $l = 3$, then ${k_1}' = k_1 - 1$, ${k_2}' = k_2 + 2$, ${k_3}' = k_3 - 1$.
Note that if $\{ k_j \}_j$ satisfies condition (\ref{cond-adpt}),
then $\{ {k_j}' \}_j$ also satisfies the condition.

Then, we need the following lemma.
\begin{lemma}\label{lem-mxp}
If $M_X$ is an even integer, then $M_{X'}$ is also an even integer. Furthermore, we have
\begin{equation*}
M_X > M_{X'}.
\end{equation*}
\end{lemma}

To prove the above lemma, we have only to consider the difference
\begin{equation*}
M_X - M_{X'} = \left( \sum_j j^2 \, k_j - \sum_j j^2 \, {k_j}' \right)  + \left( I_X - I_{X'} \right).
\end{equation*}
For each term, we have
\begin{gather*}
\sum_j j^2 \, k_j - \sum_j j^2 \, {k_j}' = 2 (l - 2) \geqslant 2,\\
I_X - I_{X'} =
\begin{cases}
0 & k_1 \leqslant n - m \; \text{or} \; k_l > m.\\
2 & \text{othrewise}.
\end{cases}
\end{gather*}
Thus, the difference $M_X - M_{X'}$ should be a positive even integer, which allows us to prove the above lemma.

If we take the sequence $X$, $X'$, $(X')'$, \ldots, then the corresponding $l$ of the sets monotonically decreases.
Thus, the sequence should end after a finite number of terms.
We denote by $X_0$ the last term of the sequence, which has at most $2$ elements.

\begin{equation*}
X_0 = \left( \left( - (l_0 - 2) - \frac{k_0}{n} \right)^{k_1^0}, \; \left( - (l_0 - 1) - \frac{k_0}{n} \right)^{k_2^0}  \right)^P
\end{equation*}
If we define $\overline{k_0} = k_0 + n (l_0 - 1)$,
then we have $- m < \overline{k_0} \leqslant n - m$,
$k_1^0 = \overline{k_0} + m$, and $k_2^0 = n - \overline{k_0} - m$, similar to (\ref{cond-adpt}).
Then, we have
\begin{equation*}
X_0 = \left( \left( 1 - \frac{\overline{k_0}}{n} \right)^{\overline{k_0} + m}, \; \left( - \frac{\overline{k_0}}{n} \right)^{n - \overline{k_0} - m}  \right)^P.
\end{equation*}
Furthermore, we have
\begin{equation}
\left( d(x, y) \right)^2 = \frac{\left( n - \overline{k_0} \right) \overline{k_0}}{n} + 2 \, i_2 \label{eq-pf_MX0}
\end{equation}

For the first assumption of Lemma \ref{lem-mx}, we have to prove that $\left( d(x, y) \right)^2$ is an even integer.
Thus, $\frac{\left( n - \overline{k_0} \right) \overline{k_0}}{n}$ should be even.
When we check every integer factor of $n$,
we can easily show that $n_0 \mid n - \overline{k_0}$ if and only if $\frac{\left( n - \overline{k_0} \right) \overline{k_0}}{n}$ is even.

Now, we take the integer $n_1$ such that $0 \leqslant n_0 n_1 < n$ and $n_0 n_1 \equiv n - \overline{k_0} \pmod{n}$.

If $n_1 = 0$, then we have $\overline{k_0} = 0$, and then $X_0 = \tilde{J}(n, m)$. By Lemma \ref{lem-mxp}, we have
\begin{equation*}
M_X > M_{X_0} = 2 m.
\end{equation*}
That is, it does not satisfy the assumption of Lemma \ref{lem-mx},
and thus we cannot add any vector of $X$ to $\tilde{J}(n, m)$ while maintaining $m$-distance.

In particular, if $n = n_0$, then by the condition $n_0 \mid n - \overline{k_0}$, we have $n_1 = 0$.
The first condition of Theorem \ref{th-main} comes from this fact.

We assume that $n_1 > 0$. Then, we have $\overline{k_0} \ne 0$, and $X_0 \ne \tilde{J}(n, m)$.
The value of $\overline{k_0}$ is determined according to three cases, as follows:
\begin{equation*}
\overline{k_0} =
\begin{cases}
n - n_0 n_1 & m < n_0 n_1.\\
n - m & m = n_0 n_1.\\
- n_0 n_1 & m > n_0 n_1.
\end{cases}
\end{equation*}
Then, we have the following forms for $X_0$:
\begin{equation*}
X_0 = 
\begin{cases}
\left( \left( \frac{n_0 n_1}{n} \right)^{n - n_0 n_1 + m}, \; \left( - 1 + \frac{n_0 n_1}{n} \right)^{n_0 n_1 - m} \right)^P & m < n_0 n_1,\\
\qquad \left( \left( \frac{m}{n} \right)^n \right) & m = n_0 n_1,\\
\left( \left( 1 + \frac{n_0 n_1}{n} \right)^{m - n_0 n_1}, \; \left( \frac{n_0 n_1}{n} \right)^{n + n_0 n_1 - m} \right)^P & m > n_0 n_1.
\end{cases}
\end{equation*}

Finally, we have only to check
\begin{equation*}
M_{X_0} = n_0 n_1 - \frac{n_1^2 n_0^2}{n} + 2 \max_{x \in \tilde{J}(n, m)} i_2 \leqslant 2 m.
\end{equation*}
Checking the above condition for the above three cases of $X_0$,
we have the following condition:
\begin{equation*}
3 n_0 n_1 - \frac{n_0^2 n_1^2}{n} \leqslant 4 m.
\end{equation*}
Furthermore, we have
\begin{equation*}
3 n_0 n_1 - \frac{n_1^2 n_0^2}{n} < 3 (n_1 + 1) n_0 - \frac{(n_1 + 1)^2 n_0^2}{n}, \qquad \text{for every} \; 0 < n_1 < n_1 + 1 < \frac{n}{n_0}.
\end{equation*}
Thus, the necessary condition for $n_1 = 1$ is that
\begin{equation*}
3 n_0 - \frac{n_0^2}{n} \leqslant 4 m.
\end{equation*}

For the other direction, if the pair $(n, m)$ satisfies the following condition,
then we can add any vector of $X_0$ for $n_1 = 1$ to $\tilde{J}(n, m)$.
This completes the proof of Theorem \ref{th-main}.
\begin{flushright}
$\square$
\end{flushright}

\section{Remarks on Theorem \ref{th-main}}\label{sec-rem}

\subsection{Case $m =4$}

We have the following lists of the sets
for which any vector of the set can be added to the corresponding $\tilde{J}(n, 4)$.
The number in brackets [ \, ] is the cardinality of the maximal four-distance sets which contain $\tilde{J}(n, 4)$.

\begin{itemize}
\item $n = 8$
\begin{itemize}
\item $\left( \left( \frac{1}{2} \right)^8 \right)$
 \qquad -- $\left( \frac{3}{2}, \left( \frac{1}{2} \right)^6, -\frac{1}{2} \right)^P$
 \qquad: $57$ vectors $[127]$.
\end{itemize}

\item $n = 9$
\begin{itemize}
\item $X^{(i)} = \left( \left( \frac{2}{3} \right)^7, \left( -\frac{1}{3} \right)^2 \right)^P$
 \qquad -- $X^{(ii)} = \left( \left( \frac{2}{3} \right)^8, -\frac{4}{3} \right)^P$
\item $X^{(iii)} = \left( \frac{4}{3}, \left( \frac{1}{3} \right)^8 \right)^P$
 \qquad -- $X^{(iv)} = \left( \left( \frac{4}{3} \right)^2, \left( \frac{1}{3} \right)^6, -\frac{2}{3} \right)^P$\\
\quad: $132$ vectors $[258]$ (conjecture).\\
\end{itemize}

\item $n = 18$
\begin{itemize}
\item $\left( \left( \frac{1}{3} \right)^{16}, \left( -\frac{2}{3} \right)^2 \right)^P$
\quad: $153$ vectors $[3213]$.
\end{itemize}

\item $n = 25$
\begin{itemize}
\item $\left( \left( \frac{1}{5} \right)^{24}, -\frac{4}{5} \right)^P$
\quad: $25$ vectors $[12675]$.
\end{itemize}
\end{itemize}

For all cases except $n = 9$, we can add all the vectors of all the sets in the list.

For the case $n = 9$, we can add all $45$ vectors of $X^{(i)}$ and $X^{(iii)}$.
However, for the sets $X^{(ii)}$ and $X^{(iv)}$, we cannot add all vectors.
We have at least $87$ vectors from $X^{(ii)}$ and $X^{(iv)}$,
which are $\left( \left( \frac{2}{3} \right)^8, -\frac{4}{3} \right) \in X^{(ii)}$ and
\begin{align*}
{X^{(iv)}}' &= \left\{ (x_1, \ldots, x_9) \in X^{(iv)} : x_i = -\frac{2}{3}, x_{j_1} = \frac{4}{3}, x_{j_2} = \frac{4}{3}, i > j_1, j_2 \right\}\\
&\quad \cup \left\{ \left( -\frac{2}{3}, \left( \frac{4}{3} \right)^2, \left( \frac{1}{3} \right)^6 \right),
\left( \frac{4}{3}, -\frac{2}{3}, \frac{4}{3}, \left( \frac{1}{3} \right)^6 \right) \right\}
\end{align*}
These are all the known vectors.

\subsection{Case $m =5$}

We have the following lists of the sets
for which any vector of the set can be added to the corresponding $\tilde{J}(n, 5)$.

\begin{itemize}
\item $n = 16$
\begin{itemize}
\item $\left( \left( \frac{1}{2} \right)^{13}, \left( -\frac{1}{2} \right)^3 \right)^P$
 \qquad: $560$ vectors $[4928]$.
\end{itemize}

\item $n = 18$
\begin{itemize}
\item $\left( \left( \frac{1}{3} \right)^{17}, -\frac{2}{3} \right)^P$
 \qquad -- $\left( \frac{4}{3}, \left( \frac{1}{3} \right)^{15}, \left( -\frac{2}{3} \right)^2 \right)^P$
 \qquad: $2466$ vectors $[11034]$.
\end{itemize}

\item $n = 25$
\begin{itemize}
\item $\left( \left( \frac{1}{5} \right)^{25} \right)^P$
 \qquad -- $\left( \frac{6}{5}, \left( \frac{1}{5} \right)^{23}, -\frac{4}{5} \right)^P$
 \qquad: $601$ vectors $[53731]$.
\end{itemize}

\item $n = 49$
\begin{itemize}
\item $\left( \left( \frac{1}{7} \right)^{47}, \left( -\frac{6}{7} \right)^2 \right)^P$
 \qquad: $1176$ vectors $[1908060]$.
\end{itemize}
\end{itemize}

For all $n$, we can add all the vectors of all the sets in the list.

\subsection{General facts}

From Theorem \ref{th-main}, we can easily show the following fact.

\noindent
{\bf Corollary \ref{th-main}.1.}
{\it For $m \geqslant 2$, the maximum $n$ such that $\tilde{J}(n, m)$ is not maximal as an $m$-distance set
is given as follows:
\begin{equation}
n = \left( 2 \left\lfloor \frac{2 (m + 1)}{3} \right\rfloor - 1 \right)^2. \label{eq-max-odd^2}
\end{equation}}

\quad

A more general version of Theorem \ref{th-main} can be shown.

\begin{proposition}\label{prop-main}
When we have $n > n_0$, then we can take an integer $n_1$ such that $0 < n_0 n_1 < n$.
If the triple $(n_0, n_1, m)$ satisfies the condition
\begin{equation}
3 n_0 n_1 - \frac{n_0^2 \, n_1^2}{n} \leqslant 4 m, \label{cond-n1},
\end{equation}
then we can add each vector of the following sets to $\tilde{J}(n, m)$ while maintaining $m$-distance:
\begin{equation}
\begin{cases}
\left( \left( \frac{n_0 n_1}{n} \right)^{n - n_0 n_1 + m}, \; \left( - 1 + \frac{n_0 n_1}{n} \right)^{n_0 n_1 - m} \right)^P & m < n_0 n_1,\\
\qquad \left( \left( \frac{m}{n} \right)^n \right) & m = n_0 n_1,\\
\left( \left( 1 + \frac{n_0 n_1}{n} \right)^{m - n_0 n_1}, \; \left( \frac{n_0 n_1}{n} \right)^{n + n_0 n_1 - m} \right)^P & m > n_0 n_1.
\end{cases}
\end{equation}
\end{proposition}

The proof of the above proposition was already discussed in the proof of Theorem \ref{th-main}.

For each example in the above proposition,
differences between distinct elements of vectors are at most $1$.
Moreover, they are all cases where differences of elements are at most $1$.
By the proof of Theorem \ref{th-main},
for every example $X$,
we have the sequence $X$, $X'$, \ldots, $X_0$.
Thus, it is easy to show that
every element of example $X$ has the form $l + \frac{n_0 n_1}{n}$ for some $l \in \mathbb{Z}$.
However, the specific form depends on $n$ and $m$.

\section{Maximal two-distance sets in $\mathbb{R}^{n - 1}$ which contain $\tilde{J}(n - 1, 2)$}\label{sec-n-1}

The Johnson graph $J(n - 1, m)$ is a subgraph of the Johnson graph $J(n, m)$.
Similar to the representation $\tilde{J}(n, m) = (1^m, 0^{n - m})^P$ of the Johnson graph $J(n, m)$ in $\mathbb{R}^{n - 1}$,
we have the following form of the representation of the Johnson graph $J(n - 1, 2)$ in $\mathbb{R}^{n - 1}$:
\begin{equation*}
\tilde{J}(n - 1, 2) = (1^2, 0^{n - 3}, 0)^{P'}.
\end{equation*}
Here, the exponent $P'$ outside the parentheses indicates that we should take every permutation of all but the last element.

If we add some vector $x = (x_1, \ldots, x_n)$ to $\tilde{J}(n - 1, 2)$ while maintaining two-distance,
then we can also add each vector of $(x_1, \ldots, x_n)^{P'}$ to $\tilde{J}(n - 1, 2)$.
The reasoning here is similar to that in the previous sections,
despite the last elements of vectors being fixed for the set.

We have the following as the form of the set of such vectors:
\begin{equation*}
X = \left( a^k, (a - 1)^{n - k - 1}, b \right)^{P'},
\end{equation*}
where we have $b = - (n - 1) a + (n - k + 1)$.

Similar to the previous sections,
for the vectors $x = (x_1, \ldots, x_n) \in \tilde{J}(n  -  1, 2)$, $y = (y_1, \ldots, y_n) \in X$,
we set $i_{j} := \left| \left\{ 1 \leqslant k \leqslant n \; ; \; x_k = 1, \; y_k = a + 1 - j \right\} \right|$ for $j = 1, 2$,
and we have 
\begin{align*}
\left( d(x, y) \right)^2 = n (n - 1) a^2 - 2 n (n - k + 1) a + (n - k + 1) (n - k + 2) + 2 \, i_2.
\end{align*}

We consider the following cases: $k = 0$, $2 \leqslant k \leqslant n - 3$, and $k = n - 2$.

For the case $2 \leqslant k \leqslant n - 3$, we have $i_2 = 0$, $1$, $2$, and so  two-distance cannot be maintained.

For the case $k = 0$, we have $i_2 = 2$, and so $\left( d(x, y) \right)^2 = 2$ and $4$.
Respectively, we have the following sets:
\begin{gather}
X_1^{\pm} = \left\{ \left( \left( \frac{2 n \pm 2 \sqrt{n}}{n (n - 1)} \right)^{n - 1}, \quad \mp \frac{2 \sqrt{n}}{n} \right) \right\}.\\
X_2^{\pm} = \left\{ \left( \left( \frac{2 n \pm \sqrt{2 n^2 + 2 n}}{n (n - 1)} \right)^{n - 1}, \quad \mp \frac{\sqrt{2 n^2 + 2 n}}{n} \right) \right\}.
\end{gather}
We have just one vector in each of the above four sets.

Similarly, for the cases $k = 1$ and $k = n - 2$, we have the following sets:
\begin{gather}
X_3^{\pm} = \left( \frac{n \pm 1}{n - 1}, \quad \left( \frac{1 \pm 1}{n - 1} \right)^{n - 1}, \quad \mp 1 \right)^{P'}.\\
X_4^{\pm} = \left( \left( \frac{3 n \pm \sqrt{10 n - n^2}}{n (n - 1)} \right)^{n - 2}, \quad \frac{- n^2 + 4 n \pm \sqrt{10 n - n^2}}{n (n - 1)}, \quad \mp \frac{\sqrt{10 n - n^2}}{n} \right)^{P'}.
\end{gather}
Here, for the sets $X_4^{\pm}$, we need $n \leqslant 10$, and we have $X_4^{+} = X_4^{-}$ for $n = 10$.
We have just $n$ vectors in each of the above four sets.

In conclusion, there are eight sets for the vectors which we can add to $\tilde{J}(n - 1, 2)$.
The distances between distinct vectors of each set are $\sqrt{2}$ or $2$;
thus, we can add all the vectors of each set to $\tilde{J}(n - 1, 2)$ while maintaining two-distance.

\begin{remark}
For each $i = 1, 2, 3, 4$, the sets $X_i^{+}$ and $X_i^{-}$ are symmetric
with respect to the hyperplane $\{ (x_1, \ldots, x_n) \in \mathbb{R}^n :  \sum_{i \leqslant n - 1} x_i = 2 \}$.
Thus, $\tilde{J}(n - 1, 2) \cup X_i^{+}$ and $\tilde{J}(n - 1, 2) \cup X_i^{-}$ are isomorphic to each other.

Furthermore, $\tilde{J}(n - 1, 2) \cup X_3^{\pm}$ is isomorphic to $\tilde{J}(n, 2)$.
$\tilde{J}(n - 1, 2) \cup X_4^{\pm}$ is just the representation of the graph in $\mathbb{R}^{n-1}$
which is made with the Seidel switching method from the Johnson graph $J(n, 2)$ and its subgraph $J(n - 1, 2)$. $($See {\upshape \cite{Se}}.$)$
In particular, only for the case $n = 5$,
the sets $J(n - 1, 2) \cup X_3^{\pm}$ and $J(n - 1, 2) \cup X_4^{\pm}$ are isomorphic to each other.
\end{remark}

Finally, we must consider all combinations of the eight sets.

We have the following list of combinations 
for which we can add all the vectors of each union to the corresponding $\tilde{J}(n - 1, 2)$.
The number in brackets [ \, ] is the total number of vectors of the union and $\tilde{J}(n - 1, 2)$.

\begin{trivlist}
\item[$(n = 5)$] \quad\\
\quad $X_1^{+} \cup X_1^{-}$:
 $\left( \left( \frac{5 \pm \sqrt{5}}{10} \right)^4, \mp \frac{2 \sqrt{5}}{5} \right)$:
$2$ vectors $[12]$.

\item[$(n = 6)$] \quad\\
\quad $X_1^{\pm} \cup X_4^{\mp}$:
 $\left( \left( \frac{6 \pm \sqrt{6}}{15} \right)^5, \mp \frac{\sqrt{6}}{3} \right)$,
 $\left( \left( \frac{9 \mp \sqrt{6}}{15} \right)^4, \frac{-6 \mp \sqrt{6}}{15}, \pm \frac{\sqrt{6}}{3} \right)^{P'}$:
$6$ vectors $[16]$.

\item[$(n = 7)$] \quad\\
\quad $X_4^{+} \cup X_4^{-}$:
 $\left( \left( \frac{21 \pm \sqrt{21}}{42} \right)^5, \frac{-21 \pm \sqrt{21}}{42}, \mp \frac{\sqrt{21}}{7} \right)^{P'}$:
$12$ vectors $[27]$.

\item[$(n = 8)$] \quad\\
\quad $X_2^{+} \cup X_4^{+}$:
 $\left( \left( \frac{1}{2} \right)^7, -\frac{3}{2} \right)$,
 $\left( \left(\frac{1}{2} \right)^6, -\frac{1}{2}, -\frac{1}{2} \right)^{P'}$:
$8$ vectors $[29]$.\\
\quad $X_2^{-} \cup X_4^{-}$:
 $\left( \left( \frac{1}{14} \right)^7, \frac{3}{2} \right)$,
 $\left( \left( \frac{5}{14} \right)^6, -\frac{9}{14}, \frac{1}{2} \right)^{P'}$:
$8$ vectors $[29]$.

\item[$(n = 9)$] \quad\\
\quad $X_1^{+} \cup X_1^{-}$:
 $\left( \left( \frac{1}{3} \right)^8, -\frac{2}{3} \right)$,
 $\left( \left( \frac{1}{6} \right)^8, \frac{2}{3} \right)$:
$2$ vectors $[28]$.\\
\quad $X_1^{+} \cup X_3^{-} \cup X_4^{-}$:
 $\left( \left( \frac{1}{3} \right)^8, -\frac{2}{3} \right)$,
 $\left( 1, 0^7, 1 \right)^{P'}$,
 $\left( \left( \frac{1}{3} \right)^7, -\frac{2}{3}, \frac{1}{3} \right)^{P'}$:
$17$ vectors $[45]$.\\
\quad $X_1^{-} \cup X_3^{+} \cup X_4^{+}$:
 $\left( \left( \frac{1}{6} \right)^8, \frac{2}{3} \right), \left( \frac{5}{4}, \left( \frac{1}{4} \right)^7, -1 \right)^{P'}$,
 $\left( \left( \frac{5}{12} \right)^7, -\frac{7}{12}, -\frac{1}{3} \right)^{P'}$:
$17$ vectors $[45]$.

\item[$(n = 17)$] \quad\\
\quad $X_1^{\pm} \cup X_2^{\mp}$:
 $\left( \left( \frac{17 \pm \sqrt{17}}{136} \right)^{16}, \mp \frac{2 \sqrt{17}}{17} \right)$,
 $\left( \left( \frac{17 \mp 3\sqrt{17}}{136} \right)^{16}, \pm \frac{6 \sqrt{17}}{17} \right)$:
$2$ vectors $[122]$.
\end{trivlist}

Note that the $29$ vectors for the case $n = 8$ and the $45$ vectors for the case $n = 9$
constitute the maximum two-distance sets introduced by Lison\v{e}k.

\quad\\

\begin{center}
{\large Acknowledgement.}
\end{center}
This paper is an outcome of the ESSP (Excellent Student in Science) project of the
Faculty of Science, Kyushu University, which was supported by the JST (Japan Science 
and Technology Agency) program Fostering Next Generation Scientists.

The first author (Bannai) gave a series of 3 hour lectures over the period Sept. 2009 to 
March 2010 (13 lectures total) to a group of six selected high school 
students. The research problems discussed in this paper were presented during the lectures. 
The second author (Sato) was one of the six members, and he solved some initial 
problems presented there in a satisfactory manner. 
The third author (Shigezumi) helped with the project on a voluntary basis by 
working on these problems independently and obtained far stronger results,
allowing the outcome of this project to be publishable as a regular research 
paper in the present form. (In Japanese, more detailed 
results written by him are available in {\tt http://researchmap.jp/j\_{}shigezumi/ESSP2009/}.) 

We sincerely thank the five other students (besides the second author): Ayaka Teshima, 
Saki Nagaki, Naoki Noda, Shuhei Noyori, and Daiki Maruyama. They also made 
considerable progress, and the discussions including them were extremely helpful. 
We thank Hiroyasu Hamada, Nobutaka Nakazono, and Yasuhiro Higashi for 
their role as teaching assistants for the class. (They wrote up the lecture notes \cite{B}.) 
We thank Etsuko Bannai for attending all the lectures and giving various 
advice. We thank Kenji Kajiwara for successfully organizing the project as 
the chairman of the mathematical part of this project in 2009.


\begin{thebibliography}{RSD}

\bibitem{B}
E. Bannai, H. Hamada, Y. Higashi, and N. Nakazono, {\it Lecture Note in ESSP 2009} (Japanese), 2010.\\
{\tt http://researchmap.jp/nakazono/ESSP2009/}

\bibitem{C}
H. T. Croft, {\it $9$-point and $7$-point configurations in $3$-space}, Proc. London Math. Soc. (3), {\bf 12} (1962),  400--424.


\bibitem{E-S2}
S. J. Einhorn and I. J. Schoenberg, {\it On euclidean sets having only two distances between points}. II, Nederl. Akad. Wetensch. Proc. Ser. A (=Indag. Math.), {\bf 69} (1966), 489-504.

\bibitem{E-F}
P. Erd\"os and P. Fishburn, {\it Maximum planar sets that determine $k$ distances}, Discrete Math., {\bf 160} (1996), 115--125. 

\bibitem{K}
L. M. Kelly, {\it Elementary problems and solutions: Solutions}: E735: {\it Isosceles $n$-points}, Amer. Math. Monthly, {\bf 54} (1947), 227--229.

\bibitem{L}
P. Lison\v{e}k, {\it New maximal two-distance sets}, J. Combin. Theory Ser. A, {\bfseries 77} (1997), 318--338.

\bibitem{Se}
J. J. Seidel, {\it Discrete non-Euclidean geometry}, in: {\it Handbook of incidence geometry}, 843--920, North-Holland, Amsterdam, 1995.

\bibitem{Sh1}
M. Shinohara, {\it Classification of three-distance sets in two dimensional Euclidean space},  European J. Combin., {\bf 25} (2004), 1039--1058.

\bibitem{Sh2}
M. Shinohara, {\it On three-distance sets in the three-dimensional Euclidean space and planar five-distance sets}, Ph.D. Thesis, Kyushu Univ., 2006.

\bibitem{Sh3}
M. Shinohara, {\it Uniqueness of maximum three-distance sets in the three-dimensional Euclidean space}, preprint.


\end{thebibliography}
\end{document}